\documentclass[oneside,11pt]{amsart}
\usepackage{amsmath}
\usepackage{amsthm}
\usepackage{quickmath}
\usepackage{thmtools}
\usepackage{pinlabel}
\usepackage{epstopdf}
\usepackage{mathtools} 
\usepackage{hyperref}
\usepackage{tikz-cd}
\usepackage{tikz}
\usepackage{amssymb}

\usepackage[T1]{fontenc}

\usepackage[small]{caption}

\theoremstyle{plain}
\newtheorem*{conjecturestar}{Conjecture}

\newcommand{\conjugation}{c}

\newcommand{\hQn}{\hat{\mathcal{Q}}_N}

\newcommand{\hQ}{\hat{\mathcal{Q}}}
\renewcommand{\-}{^{-1}}
\DeclareMathOperator{\charc}{char}

\newcommand{\ios}{\iota(\sigma)}

\date{}

\begin{document}

\title{Acylindrical Hyperbolicity of Out($W_n$)}

\author{Brendan Burns Healy}
\address{Department of Mathematical Sciences\\
         University of Wisconsin--Milwaukee\\
         PO Box 413\\
         Milwaukee, WI 53211\\
	 USA}
\email{healyb@uwm.edu}

\begin{abstract}
We prove that the group of outer automorphisms of the free Coxeter group $W_n$ is acylindrically hyperbolic in the sense of Osin. As an application, we observe that any $\CAT(0)$ space admitting a geometric action by $\Out(W_n)$ must contain a rank-one geodesic. The theorem proceeds from expanding on a well-known relationship between $\Out(W_n)$ and the outer automorphism group of free groups.
\end{abstract}

\date{\today}

\maketitle

\begin{section}{Introduction}

Coxeter groups, as abstractions of reflection groups, play a critical role in both classical geometry and Lie theory, and the universal example, the free Coxeter group $W_n := \ast_n \Z / 2\Z$, has long been studied by geometric group theorists, dating back more than 20 years to \cite{muhl}, for its importance to \emph{Coxeter systems}. The free Coxeter group contains an easy-to-see characteristic copy of $F_{n-1}$, which induces a close relationship between the outer automorphisms of these groups, which we detail further in Section 3. The group $\Out(W_n)$ has been studied by, among others, Gutierrez, Piggott, and Ruane in \cite{gprr}, and is expected to inherit many of the properties enjoyed by $\Out(F_{n-1})$, though with some key differences. \\
\indent The notion of \emph{acylindrical hyperbolicity} was introduced by Osin in 2013 to unify various ideas introduced previously, such as weak proper discontinuity and weakly contracting elements. In \cite{osin}, a short survey is provided of known examples of groups which satisfy this property, including most mapping class groups and 3-manifold groups which are not virtually polycyclic, as well as $\Out(F_n)$. Because of the relationship between this last example to $\Out(W_n)$, we might expect this group to enjoy this property. Our main theorem is that this is true.

\begin{restatable*}{theorem}{outah}
\label{theorem:outah}
$\Out(W_n)$ is acylindrically hyperbolic, for $n \geq 3$.

\end{restatable*}

The property of acylindrical hyperbolicity for a group is of interest for many reasons, some of which are explored in \cite{osin}. One particularly strong consequence is that these groups are \emph{SQ-universal}, meaning that any countable group embeds in some quotient of the group. Furthermore, acylindrical hyperbolicity is a generalization of relative hyperbolicity. However it is known that $\Out(W_n)$ is not relatively hyperbolic - indeed it is thick due to \cite{das}. As such, we do not get this result as a consequence of known properties for this group. Furthermore this is behavior that actually distinguishes our group from outer automorphism groups associated to one-ended hyperbolic groups. Levitt finds in \cite{levitt} that $\Out(G)$, for $G$ one-ended hyperbolic and not a surface group, virtually admits a short exact sequence where the normal subgroup is free abelian. As acylindrical hyperbolicity passes to infinite normal subgroups by \cite{osin}, this is an obstruction to $\Out(G)$ being acylindrically hyperbolic.

Gersten in \cite{gersten} demonstrates that $\Out(F_n)$ is not $\CAT(0)$ for $n \geq 4$ and Bridson and Vogtmann proved that $\Out(F_n)$  has exponential Dehn function~\cite{bridvogt}.
The particular automorphisms used in both of these proofs involve transvections, which are not in the image of the map $\iota\colon \Aut(W_n) \to \Aut(F_{n-1})$ 
defined in the next section, 
%
%
leaving the possibility open that $\Out(W_n)$ is a CAT(0) group.
However, piecing together known results of Charney--Sultan and Sisto, 
we are able to derive strong restrictions on the class of CAT(0) spaces 
which could potentially admit a geometric $\Out(W_n)$-action. \\

\noindent {\bf Corollary \ref{outcat} } {\it Suppose $\Out(W_n)$ acts geometrically on $X$ a CAT(0) space for $n \geq 3$. Then $X$ contains a rank-one geodesic. In particular, $\Out(W_n)$ cannot act geometrically on a Euclidean building. } \\

The author would like to thank Jeffrey Carlson for helpful comments and suggestions.

\end{section}

\begin{section}{ \texorpdfstring{$\Out(W_n)$}{OutWn} }
\label{sec:OW}
In this section we discuss the relationship between $\Out(W_n)$
and $\Out(F_{n-1})$. The results from here up to and including \ref{inclusion} are known and in the literature,
but we provide some proofs for convenience.

Throughout this work we
will denote by $W_n = \langle w_i \mid w_i^2 \rangle$
the free Coxeter group on $n$ letters
and by $G$ the subgroup of words of even length,
which is generated by the elements $x_i := w_0 w_i$
for $1 \leq i \leq n-1$.
We will momentarily show $G$ is free on the $x_i$ and so
we will often abusively denote this subgroup by $F_{n-1}$.

\begin{definition}
Let $\phi \in \Aut(W_n)$. Then we call $\phi$ a \emph{partial conjugation} if for some proper nonempty subset $J \subset \{ 0, 1, \ldots n-1\}$ and $z \in W_n$
 
$$
\phi(w_i) = \left\{
   \begin{array}{ll}
      zw_iz^{-1} & \quad j \in J \\
      w_i & \quad \text{else}
    \end{array}
  \right.
$$

We will often refer to the automorphism above as $F^z_J$. An \emph{elementary partial conjugation} is one where the cardinality of $J$ is 1 and $z=w_i$ where $\{i\} \neq J$. Any partial conjugation is the product of elementary partial conjugations.
\end{definition}

We will make repeated use of the following result, described originally in the general case by \cite{collins} and written in the following form for $W_n$ in \cite{gprr}.

\begin{theorem} [\cite{collins}] 
\label{thm:decomposition}
There is a decomposition
\[ \Aut(W_n) = \Aut^0(W_n) \rtimes \Sigma_n \]
where $\Aut^0$ is the group generated by partial conjugations and $\Sigma_n$ is the full symmetric group on $n$ letters acting by permutation of the generators $w_i$. Moreover, we get the following split exact sequence, where the copy of $W_n$ is the subgroup of inner automorphisms.
\[
	1 \to W_n \to \Aut^0(W_n) \to \Out^0(W_n) \to 1
\]
\end{theorem}

Note that in the above short exact sequence, we are \emph{defining} $\Out^0(W_n)$ to be the group that fits into that sequence. We prove the following here for completeness, due to its importance in the rest of the article.

\lem[{\cite{muhl}}]\label{thm:free}
Consider the subgroup $\langle w_0 w_i \mid i \in \{1, \ldots n-1\} \rangle$ of $W_n$. 
This group is free on $n-1$ generators.
\elem

\pf
It is clear that none of the $n-1$ generators are redundant. So we are reduced to demonstrating there are no relations. We begin by noting that \[(w_0 w_i)^{-1} = w_i w_0 \]

Cancellation can only happen in the form $w_i w_i = 1$, as these are the only relators in $W_n$. This is only the case if we have $w_i w_0 w_0 w_j$, in which case this is equal to $w_i w_j$, which is irreducible and simply another expression for the element $x_i x_j$, or $w_0 w_j w_j w_0$, which translates to $x_{j-1} x_{j-1}^{-1}$, i.e., a free reduction in $F_{n-1}$.
\epf

Another way of recognizing the above subgroup is exactly the subgroup of even length words in the given generating set.

\begin{lemma}[{\cite{pigruane}}]\label{thm:characteristic}
	The subgroup in Lemma~\ref{thm:free} is characteristic. 
\end{lemma}
\begin{proof}
	Since $\Aut(W_n)$ is generated by permutations of generators 
	and partial	conjugations by \ref{thm:decomposition}, 
	all automorphisms preserve parity of word length
	and hence $G$.
\end{proof}

This lemma allows us to restrict automorphisms to obtain a map
\begin{equation}\label{def:iota}
\iota\colon \Aut(W_n) \rightarrow \Aut(F_{n-1})
\end{equation}

\lem [\cite{muhl}]
\label{inclusion}
This map is injective for $n \geq 3$.
\elem

\pf
Let $\phi$ lie in the kernel of $\iota$, meaning it fixes each $w_0 w_i$. We show $\phi$ fixes all $w_i$. Indeed, it is enough to show it fixes $w_0$; as then it clearly fixes $w_i = w_0 \cdot w_0 w_i$ as well.

Note each $\phi(w_i)$ must start and end with the same letter,
because $w_i$ does and 
$\phi$ is a composition of partial conjugations and permutations
by \ref{thm:decomposition}.
%
Write $z$ for a reduced word representing $\phi(w_0)$.
By assumption, $\phi(w_0w_i) = w_0 w_i$,
so $\phi(w_i) = z^{-1}w_0 w_i$. \\
\indent If this word is reduced, it ends with $w_i$, thus it must also begin with $w_i$ meaning $z$ ends with $w_i$.
Now, if $z$ did not begin with $w_0$,
the word $z^{-1}w_0 w_i$ would already be reduced
and hence, impossibly, $z$ would end with $w_i$ for each $i$,
so $z = \phi(w_0)$ begins and ends with $w_0$.
It follows $z^{-1} w_0$ is either $1$
or also begins with $w_0$.
But $z^{-1}w_0 w_i$ begins with $w_i$,
so indeed $z^{-1} w_0 = 1$, or $\phi(w_0) = w_0$. \\
\indent If the word $z^{-1}w_0 w_i$ is not reduced, then $z$ begins with $w_0$. We may also assume it starts with $w_0 w_i$, or else we would be in the above case as it still ends (and begins) in $w_i$. But again, it must impossibly start with $w_0 w_i$ for all $i$ simultaneously. \qedhere

\epf

Because our goal is to say something about $\Out(W_n)$, we look at what happens to elements of $\Inn(W_n)$. While it is not quite true that they map into $\Out(F_{n-1})$, we find this is \emph{almost} the case.

\lem
Let $r \in \Aut(F_n)$ be defined by $r(x_i) = x_i^{-1}$ for all $i$. Then
\[ \iota (\Inn(W_n) ) \subset \Inn(F_{n-1}) \rtimes \{r\}. \]
\elem

\pf
For any element $w \in W$ write $\conjugation_w$ for conjugation by that element.
We must show the restriction of each $\conjugation_{w_i}$ 
to $F_{n-1}$
can be written in terms of $r$ and $c_{x_i}$.
%
\begin{itemize}
\item First, $\conjugation_{w_0} (w_0 w_j) = w_j w_0 = (w_0 w_j)^{-1}$ for all $j$,
so $\iota(\conjugation_{w_0}) = r$.
\item For $j > 0$, note that since $x_j = c_{w_0}\big(r(x_j)\big)$ we have
\[
	\conjugation_{w_i} (x_j) = 
	\conjugation_{w_i} \conjugation_{w_0}r(x_j) =
  \conjugation_{x_i}^{-1} r(x_j).
  \qedhere
\]
\end{itemize}
\epf

Thus $\iota$ descends to a well-defined map
$\bar\iota\colon \Out(W_n) \to \Out(F_{n-1})/\langle \langle r \rangle \rangle$,
where $\langle \langle r \rangle \rangle$ denotes the normal closure of 
$\langle r \rangle$. This relationship is summarized in the following commutative diagram.
\begin{figure}[h]
\begin{center}
\begin{tikzcd}
\Aut(W_n) \arrow{r}{\iota} \arrow{dd}{\bar{q}}
& \Aut(F_{n-1}) \arrow{d}{q}\\
 & \Out(F_{n-1}) \arrow{d}{q_r} \\
\Out(W_n) \arrow[dotted]{r}{\bar{\iota}}
& \Out(F_{n-1})/\langle \langle r \rangle \rangle 
\end{tikzcd}
\end{center}
\caption[Relating Outer Automorphisms]{Diagramatic relations of the groups}
\label{diagram}
\end{figure}

\noindent We will show in fact that the kernel of
the restriction of $q_r$ to $q_r^{-1}(\im \bar\iota)$ 
is the subgroup of order two generated by $r$.

\lem
\label{thm:comm} 
There is a containment
$[r, \im \iota] \subset \Inn(F_{n-1})$.
\elem

\pf 
It will be enough to show the containment 
on a set of generators of $\im \iota$ because of the general identity $[r,ab] = [r,a]\cdot a[r,b]a\-$ and because
$\Inn(F_{n-1})$ is normal in $\Aut(F_{n-1})$.
By \ref{thm:decomposition}, elements of $\Sigma_n$ and elementary partial conjugations will suffice. We proceed by cases.

\begin{itemize}
\item The symmetric group is generated by transpositions $\sigma = (w_i w_j)$.
\begin{itemize}
\item $\sigma = (w_i w_j), i \neq 0 \neq j$. Then it is easy to see $\iota(\sigma)$ permutes $x_i, x_j$ and that this map commutes with inverting every generator.
\item $\sigma = (w_0 w_i)$. In this case, $\sigma(w_0 w_ j) = w_i w_j = w_i w_0 w_0 w_j$ for $i \neq j$ and $\sigma(w_0 w_i) = w_i w_0$. Then $\ios(x_j) = x_i^{-1} x_j$ for $j \neq i$ and $\ios(x_i) = x_i^{-1}$. This action commutes with $r$ as well.
\end{itemize}
\item $\phi = \phi^{-1}$ is a partial conjugation which takes $w_i \mapsto w_k w_i w_k$ and fixes other generators.
\begin{itemize}
\item $k \neq 0 \neq i$. Then $\iop$ fixes all generators of $F_{n-1}$ except $x_i$, and $\iop(x_i) = x_k x_i^{-1} x_k$.  We see that $[r,\iop] (x_j) = x_j$ for all $j$.

\item $k=0$. In this case, $\iop$ inverts $x_i$ and fixes the other free generators. This automorphism commutes with inverting all generators.
\item $i=0$, so $\phi(w_0) = w_k w_0 w_k$. Quick calculations show that $\iop (x_j) = x_k^{-2} x_j$. Then one can check that $(\iop \circ r \circ \iop \circ r ) (x_j) = x_j^{x_k^{-2}}$  \qedhere

\end{itemize}
\end{itemize}
\epf

Thus the class $q(r)$ of $r$ in $\Out(F_{n-1})$
lies in the center of the subgroup $\im(q \circ \iota)$
and particularly is normal,
so that
%
$\langle r \rangle = \langle \langle r \rangle \rangle \cap \im(\iota \circ q)$.
More to the point, this allows us to replace Figure \ref{diagram} with Figure \ref{diagram2}.

\begin{figure}[h]
\begin{center}
\begin{tikzcd}
\Aut(W_n) \arrow{r}{\iota}[swap]{\cong} \arrow{dd}{\bar{q}}
& \im(\iota) \arrow[hook]{r} \arrow{d} 
& \Aut(F_{n-1}) \arrow{d}{q} \\ 
& \im(\iota \circ q) \arrow{d} \arrow[hook]{r} 
& \Out(F_{n-1}) \arrow{d}{q_r} \\
\Out(W_n) \arrow[dotted]{r}{\bar{\iota}}[swap]{\cong}
& \im(\iota \circ q) / \langle r \rangle \arrow[hook]{r}& \Out(F_{n-1})/\langle \langle r \rangle \rangle 
\end{tikzcd}
\end{center}
\caption[A modified relationship diagram]{Involution normality in the image}
\label{diagram2}
\end{figure}

\end{section}

\newpage

\begin{section}{Relationship to \texorpdfstring{$\Out(F_n)$}{OutFn} }
\label{sec:Irred}

The groups $\Out(F_n)$ are very widely studied in part due to their analogous behavior to mapping class groups. One important guiding question in this investigation is what space should play the role corresponding to the curve complex. To understand one relevant answer to this question, we review some terms from the study of these groups.


\begin{definition}
 A \emph{free factor} of $F_n$ is a proper subgroup $A < F_n$ such that there exists some other subgroup $B < F_n$ where $F_n = A \ast B$.
\end{definition}

\begin{definition}
 An outer automorphism $\psi$ of $F_n$ is called \emph{fully irreducible} if there is no free factor $A$ such that some power $\psi^p$ fixes the conjugacy class of $A$.
\end{definition}

In \cite{bestfei} Bestvina and Feighn demonstrate that a natural complex, called the \emph{free factor complex} whose vertices correspond to conjugacy classes of free factors of $F_n$, admits an action by $\Out(F_n)$ and is $\delta$--hyperbolic. While this action is not geometric, it is demonstrated that fully irreducible automorphisms act with weak proper discontinuity (WPD) - a property with an intimate relationship to acylindricity. In fact, by \cite{osin}, this is sufficient to find that $\Out(F_n)$ is acylindrically hyperbolic and that these elements are generalized loxodromics.

In order to eventually demonstrate non-elementarity of our action, we will need to find a fully irreducible element in the image of our map $\bar{\iota}$. Because these elements will actually live in a quotient of the group, we will prove in Lemma \ref{lemma:NormalSubgroup} that the images of generalized loxodromics under finite kernel quotients will be generalized loxodromics for the quotient group. We follow the lead of Gersten and Stallings to find such a fully irreducible automorphism.

\begin{lemma}[\cite{gerstal}]
\label{lem:Gersten}
Let $\phi$ be an outer automorphism of $F_n$ and let $\bar{\phi}$ be the associated matrix which transforms the standard basis elements of $\mathbb{Z}^n = F_n^{ab}$. Then if the characteristic polynomial for $M_\phi$ is irreducible over $\mathbb{Q}$ and the matrix itself is primitive (meaning it has non-negative entries and some positive power has all positive entries), the automorphism is fully irreducible.
\end{lemma}

The way we use this lemma is inspired by examples stated later in \cite{gerstal}. 
In particular an automorphism with the following $(n-1) \times (n-1)$ matrix meet the above criteria:
\[
A = \begin{bmatrix} 
  0 & 1 & 0 & \cdots & 0\\
  0 & 0 & 1 & \cdots& 0 \\
  \vdots &\vdots& \vdots &\ddots& \vdots \\
  0 & 0 & 0 & \cdots & 1 \\
  1 & 4 & 4 & \cdots & 4  
  \end{bmatrix},
\]
The $(n-1)$ power of these matrices will have all positive entries
and elementary matrix operations show 
\begin{eqnarray*} 
\charc(A) &=& x^{n-1} - 4x^{n-2} - 4x^{n-3} - \cdots - 4x -1,
\end{eqnarray*}
which has no rational roots.

\begin{lemma}
\label{iwip}
The image of the map
\[ \bar{\iota} : \Out(W_n) \hookrightarrow \Out(F_{n-1}) / \langle \langle r \rangle \rangle \]
contains an element which factors through an outer automorphism that is fully irreducible.
\end{lemma}

\pf
We will define an automorphism whose induced matrix is as above.
From the proof of \ref{thm:comm}
we see there is an element of $\im \iota \subset \Aut (F_{n-1})$
inverting any one $x_i$ and fixing the other $x_j$,
and another element taking $x_i \mapsto x_k x_i\- x_k$ and fixing the other $x_j$.
Composing these yields an element $c_{i,k}$ of $\im \iota$
sending $x_i \mapsto x_k x_i x_k$ and fixing the other $x_j$.

The subgroup of $\Sigma_n < \Aut(W_n)$
fixing $w_0$ induces an action of the subgroup 
$\Sigma_{n-1}$ of $\Aut(F_{n-1})$ permuting the 
generators $x_i$. Then the matrix corresponding to the compositions 
\begin{eqnarray*} 
\phi_A &=& c_{1,n-1 }^2\circ \cdots c_{1,3}^2 \circ c_{1,2}^2 \circ (1 \ 2 \ \cdots \ n-1),
\end{eqnarray*}
is $A$ as above.

\epf

\end{section}


\begin{section}{Acylindrical Hyperbolicity}

We review what it means for a group action to be acylindrical, in order to show this for our constructed action.

\defi
 An metric space action $G \curvearrowright S$ is called \emph{acylindrical} if for every $\epsilon > 0$ there exist $R(\epsilon ),N(\epsilon ) > 0$ such that for any two points $x,y \in S$ such that $d(x,y) \geq R$, the set
$$ \{ g \in G \mid d(x,g.x) \leq \epsilon, d(y,g.y) \leq \epsilon \}$$
has cardinality less than $N$.
\edefi

Recall a group is \emph{acylindrically hyperbolic} if it acts acylindrically and non-elementarily on a hyperbolic space. 

\defi
Let $G$ be an acylindrically hyperbolic group.  An element $g \in G$ is called a \emph{generalized loxodromic} if there's an acylindrical action $G \curvearrowright X$ for $X$ hyperbolic such that $g$ acts as a loxodromic.  An individual isometry $g$ is a \emph{loxodromic} if for some basepoint $x_0 \in X$, the map $\Z \rightarrow X$ defined by $n \mapsto g^n . x_0$ is a quasi-isometry.
\edefi

It should be noted that part of the following lemma is proven in \cite{minosin} as Lemma 3.9. We prove a version here in order to emphasize the relationship between the generalized loxodromics.

\begin{lemma}
\label{lemma:NormalSubgroup}
Let $G$ be a group with a finite normal subgroup $N \vartriangleleft G$. Then $G/N$ is acylindrically hyperbolic if and only if $G$ is acylindrically hyperbolic. Furthermore the generalized loxodromics of $G/N$ are exactly the images of the generalized loxodromics of $G$ under the quotient map.
\end{lemma}

\begin{proof}
For the first direction, suppose we know $G$ is acylindrically hyperbolic and $\gamma \in G$ is a generalized loxodromic. By the characterization of acylindrical hyperbolicity [Theorem 1.2,AH1] in \cite{osin} there is a hyperbolic graph that the group acts on acylindrically and non-elementarily with bounded quotient, such that $\gamma$ acts loxodromically. Call such a graph $\mathcal{Q}$.

The first thing we will do is make a slight modification to $\mathcal{Q}$. Unlike in uniquely geodesic spaces such as $\CAT(0)$ spaces, fixed point sets in arbitrary $\delta-$hyperbolic spaces are not as nice as we would like, so we will add in a little extra structure. Let $\delta$ represent the constant of hyperbolicity for $\mathcal{Q}$ and define $\tau := \max \{ d(x,n.x) | x \in X^{(0)}, n \in N \}$. This second value is guaranteed to be finite because there is only one vertex in the bounded quotient, and $N$ is a finite, normal subgroup.

Define $\hat{\mathcal{Q}} := \mathcal{Q} \cup E \cup F$. The set $E$ consists of new combinatorial edges of length $ \tau |N|$ between any two vertices in $\hQ$ of the form $v, n.v$ for $n \in N$, and the set $F$ will consist of copies of an $|N|$--simplex metrized in the usual way such that the 1-skeleton of the boundary is the $N$--orbit of the set of new combinatorial edges emanating from a given point. We note that these two spaces are quasi-isometric by observing that $\mathcal{Q}$ embeds into $\hat{\mathcal{Q}}$ in the natural way such that distances are not changed (taking any path through these new simplices is not a shortcut, by design), and the embedding is quasi-onto with a constant no larger than the diameter of our new, finite-dimensional simplices. The group $G$ will act on $\hat{\mathcal{Q}}$ in the natural way on the embedded copy $\mathcal{Q}$ and permute the edges in $E$ according to their endpoints. Finally, for the higher dimensional simplices, extend the action on their boundaries to the interior isometrically. 

We will label the fixed point set of $N$ under this action $\hQn$, and observe that the inclusion $\hQn \hookrightarrow \mathcal{Q}$ is a $G$--equivariant map. As $N$ is in the kernel of the action on $\hQn$ by design, this action factors through $G/N \curvearrowright \hQn$. In order to show that this action on $\hQn$ demonstrates acylindrical hyperbolicity for the group, we must show three things:
\begin{enumerate}
\item $\hat{\mathcal{Q}}_N$ is hyperbolic.
\item This action is non-elementary.
\item This action satisfies acylindricity.
\end{enumerate}

For the first task we recall that hyperbolicity is a quasi-isometry invariant, so we know that $\hQn$ is hyperbolic. Because of the metric we equipped it with, the natural injection is a quasi-isometric embedding. To get quasi-density we observe that every element of $N$ fixes each $|N|$-simplex set-wise. By the way we defined the action, we may use the Brouwer fixed point theorem to observe that there is some point in each such simplex that is in the subspace $\hQn$. Therefore, this embedding is also quasi-dense by the assumption that every point was in the 0-skeleton of some added $|N|$-simplex.

To demonstrate non-elementarity, observe that if an element $g \in G$ acted loxodromically on $\mathcal{Q}$, it also acted loxodromically on $\hQ$. This implies that the map, for any $x_0 \in \hQ$, $\Z \hookrightarrow \hQ$ defined by $n \mapsto g^n.x_0$ was a quasi-isometric embedding. If we compose this map with the quasi-inverse to the embedding of $\hQn \subset \hQ$ (which we showed was a full quasi-isometry), we get a quasi-isometric embedding $\langle [g] \rangle \cong \Z \hookrightarrow \hQn$, demonstrating that $[g]$ must act loxodromically, where $[g]$ is the class of $g$ in the quotient map $G \rightarrow G/N$. Therefore, the set of elements which act loxodromically is preserved, ensuring non-elementarity. We also note that the image of the element $\gamma$ will be a loxodromic under this action.

For acylindricity we begin by letting $R(\epsilon), N(\epsilon)$ be constants depending on $\epsilon$ that demonstrate the acylindricity of the action $G \curvearrowright \hQ$. By quasi-isometry, points at distance $\ell$ in $\hQn$ were at distance no more than $c \ell + c$ in $\hQ$ where $c$ is the quasi-isometry constant. So for new constants $R', N'$, we can take $R'(\epsilon) = R(c \epsilon + c)$, and $N'(\epsilon) = N(\epsilon)$, noting that the class of generalized loxodromics was preserved under the quotient, and that the map $G \rightarrow G /N$ was finite-to-one and onto, implying there are no \emph{more} elements fixing nearby points.

Now suppose we have that $G/N$ is acylindrically hyperbolic with some generalized loxodromic $\gamma$. Again, let $\mathcal{Q}$ be such a graph with a cobounded, acylindrical action by $G/N$. We can define an action on the same space by $G$ by allowing all points to have stabilizer $N$, which is normal. We have not changed the space, so $\mathcal{Q}$ is still hyperbolic.

The action is non-elementary by the observation that if $gN$ represents a generalized loxodromic element of $G/N$, then any pre-image $g$ acts identically on $\mathcal{Q}$. Again, this tells us that the pre-image of generalized loxodromics are generalized loxodromics in the new action, and in particular will be true for our chosen $\gamma$. 

Acylindricity is obtained by taking the same $R(\epsilon)$ for the original, but now allowing the size of group elements to be multiplied by the order of $N$, to count all pre-images of elements which will have the same action on $\mathcal{Q}$. This means that $N'(\epsilon) = |N| N(\epsilon)$, which was finite by assumption.

Because we showed that any action of either group where a given generalized loxodromic acts loxodromically generates an action of the other where its (pre-)image does so, we get the desired correspondence of generalized loxodromic elements. Though a universal action, one in which all generalized loxodromics act as loxodromics, is not guaranteed to exist by \cite{abbott}, for any generalized loxodromic element there is some hyperbolic space for which it acts loxodromically. We can then use this argument to get a space with the correct action for the quotient/extension.
\end{proof}

\section{Acylindrical Hyperbolicity of Outer Automorphism Groups of Virtually Free Groups}

Our main theorem of the paper is the following.

\outah

\begin{proof}
The group $\Out(W_3)$ is virtually a non-abelian free group, from \cite{collins}. As such it is hyperbolic and non-elementary, so it is acylindrically hyperbolic.

For $n \geq 4$ Lemma \ref{iwip} tells us that there is an element in the image of $\iota$ that descends to a fully irreducible automorphism. From the relationship summarized in Figure \ref{diagram2}, the group $\im(\iota \circ q)$ is a subgroup of $\Out(F_{n-1})$, which is acylindrically hyperbolic for $n \geq 4$. Furthermore, it is a subgroup which contains a generalized loxodromic element and is not virtually cyclic, so it too is acylindrically hyperbolic. Any element which is fully irreducible will be a generalized loxodromic, by \cite{bestfei} and \cite{osin}. We can then apply Lemma \ref{lemma:NormalSubgroup} to assert that the quotient $\im(\iota \circ q) / \langle r \rangle$ is acylindrically hyperbolic where generalized loxodromics include the images of fully irreducible automorphisms. Because $\Out(W_n)$ is isomorphic to its image under this map, this proves the theorem.
\end{proof}

\begin{corollary} \label{outcat} Suppose $\Out(W_n)$ acts geometrically on $X$ a CAT(0) space. Then $X$ contains a rank-one geodesic. In particular, $\Out(W_n)$ cannot act geometrically on a Euclidean building or symmetric space.
\end{corollary}
\begin{proof}

This follows from Theorem 4.1 and the fact stated in the appendix. 

\end{proof}

It seems at first that a similar result should hold for arbitrary finite free products of finite groups. Collins in \cite{collins} constructs a characteristic free subgroup of $G = \ast_{i=1}^n A_i$ for $A_i$ finite, which arises as the kernel of the canonical projection $G \twoheadrightarrow \oplus A_i$. Indeed he goes on to prove that this subgroup induces an injection (in the case of sufficient complexity) of $\Aut(G) \hookrightarrow \Aut(F_k)$. However, this free group will have rank considerably larger than $n$ and the image of this map fails to contain automorphisms satisfying the hypotheses of Lemma~\ref{lem:Gersten}. In the case of $W_n$, the subgroup induced in this way is distinct from that created in Section 2.

However, this obstruction does not guarantee that such groups are decidedly not acylindrically hyperbolic. In fact one might suspect that all virtually free groups have acylindrically hyperbolic outer automorphism group. We see that this is not the case. Pettet in \cite{Pettet} studies a class of groups coming from graphs of groups with finite-by-cyclic vertex groups and finite edge groups such that the resulting group is virtually free. In particular he demonstrates the following, which precludes the associate $\Out(G)$ from being acylindrically hyperbolic by the non-elementary condition.

\begin{theorem}[\cite{Pettet}]
There exist virtually free groups whose outer automorphism groups are finite.
\end{theorem}

We can combine this with an easy to state condition from \cite{boutin} which ensures the above property. Therefore, we get the following corollary.

\begin{corollary}[\cite{Pettet},\cite{boutin}]
Let $G$ be a group given by the following exact sequence, with $K$ finite: $1 \rightarrow F_n \rightarrow G \rightarrow K \rightarrow 1$. The conjugation action of $G \rightarrow F_n$ induces a homomorphism $\theta: K \rightarrow \Out(F_n)$. If the centralizer of $\theta(K)$ in $\Out(F_n)$ is finite, then $\Out(G)$ is not acylindrically hyperbolic.
\end{corollary}

In fact we can find non-examples solely among the class of right-angled Coxeter groups.

\begin{proposition}[\cite{gprr},\cite{virtfree}]
Let $\Gamma$ be a finite graph which is homeomorphic to $[0,1]$. Then $W_\Gamma$ is virtually free and $\Out(W_\Gamma)$ is finite.
\end{proposition}

\begin{proof}
The fact that $W_\Gamma$ is virtually free follows from the achordal condition on $\Gamma$ as in \cite{virtfree}. Furthermore, because this graph cannot have any separating intersections of links (SILs) as defined in \cite{gprr}, we get that it has finite outer automorphism group.
\end{proof}

A larger class of non-examples exists. A generalization of the SIL condition was defined by Sale and Susse in \cite{ss}. In that paper they characterize when the outer automorphism group will be virtually abelian, and therefore not acylindrically hyperbolic. Again by combining this with the main result from \cite{virtfree}, we can conclude the following.

\begin{proposition}[\cite{ss},\cite{virtfree}]
Let $\Gamma$ be a finite graph which is achordal and contains no FSILs or STILs as in \cite{ss}. Then $W_\Gamma$ is virtually free and $\Out(W_\Gamma)$ is not acylindrically hyperbolic.
\end{proposition}

\end{section}

\appendix

\begin{section}{ \texorpdfstring{$\CAT(0)$}{CAT0} vs. Acylindrically Hyperbolic}
Here we make an observation following from existing work, that we find useful to note.

\fct [Charney-Sultan, Sisto]
Let $G$ be a group, which is not virtually cyclic, acting geometrically on a CAT(0) space $X$. Then $G$ is acylindrically hyperbolic if and only if it contains an element $g$ which acts as a rank-one isometry on $X$. Furthermore, the set of generalized loxodromics is precisely the set of rank-one elements.
\label{equiv}
\efct

\pf
($\Leftarrow$)
This follows from Theorem 5.4 in \cite{bestfuji}, where it is stated that a geodesic in a CAT(0) space is contracting exactly when it fails to bound a half-flat, meaning rank-one geodesics are contracting. Next, contracting elements are also \emph{weakly contracting} in the sense of \cite{sisto}. Sisto goes on to prove ~\cite[Theorem~1.6]{sisto} that such an element is contained in a virtually cyclic subgroup, labelled $E(g)$, which is hyperbolically embedded in the group. This is one of four equivalent conditions for being a generalized loxodromic, as per \cite[Theorem~1.4]{osin}.

($\Rightarrow$) If $G$ is acylindrically hyperbolic, then it contains at least one generalized loxodromic. 
Indeed, by assumption $G$ admits an action on a hyperbolic space 
under which no element acts parabolically, so by non-elementarity of the action, at least one element must act loxodromically. 
Call this element $g$. We know by a result of Sisto that $g$ is Morse in $G$ \cite{sisto2}. 
But an equivalence in the setting of CAT(0) groups proved in \cite{charney} states that a (quasi-)geodesic in a CAT(0) space is contracting if and only if it is Morse and if and only if it is of rank one. 
As our action is geometric, as our element $G$ is Morse,
its axes are as well, so $g$ acts as a rank-one isometry. 
\epf

One reason for noting this equivalence is it allows us to recast traditionally geometric CAT(0) statements in terms of acylindrical hyperbolicity. For example, using \ref{equiv} one can restate of Ballman and Buyalo's rank rigidity conjecture as follows:

\begin{conjecturestar}{\emph{Rank Rigidity}} \cite{BB} \\
Let $X$ be a locally compact geodesically complete CAT(0) space and $G$ a discrete group acting geometrically on $X$. If $X$ is irreducible, then either
\begin{itemize}
\item $X$ is a Euclidean building or higher rank symmetric space 
\newline or
\item $G$ is acylindrically hyperbolic.
\end{itemize}
\end{conjecturestar}

\end{section}

\bibliography{refs}
\bibliographystyle{alpha}
\end{document}